\documentclass[onefignum,onetabnum,noabbrev]{siamonline171218}

\usepackage{amsmath}
\usepackage{amssymb}
\usepackage[textsize=tiny]{todonotes}
\usepackage{siunitx}
\usepackage{algpseudocodex}
\usepackage{graphicx}
\usepackage[margin=2.5cm]{geometry}
\usepackage{hyperref}
\usepackage{tikz}
\usepackage{drawmatrix}
\usepackage{xcolor}
\usetikzlibrary{positioning}
\newcommand{\Ex}{\mathbb{E}}

\newcommand{\G}{\mathbf{G}}
\newcommand{\A}{\mathbf{A}}
\newcommand{\precmat}[1]{P_{#1}}
\usepackage{subcaption}

\definecolor{nottblue}{HTML}{384257}

\definecolor{lightgray}{HTML}{f2f9ff}

\definecolor{lightorange}{HTML}{26bb9c}
\definecolor{diagram}{HTML}{00939f}

\newcommand\blfootnote[1]{%
  \begingroup
  \renewcommand\thefootnote{}\footnote{#1}%
  \addtocounter{footnote}{-1}%
  \endgroup
}

\author{Victor Trappler\footnote{Inria, CNRS, Univ. Grenoble-Alpes, Grenoble-INP, LJK, 38000 Grenoble, France}  \footnote{AI4Sim, Eviden BDS R\&D, Echirolles, France}  \footnote{Current Affiliation: \'Ecole Centrale de Lyon, CNRS UMR 5208, Institut Camille Jordan, 36
Avenue Guy de Collongue, 69134 \'Ecully, France} \and Arthur Vidard\footnotemark[1]}

\title{State-dependent preconditioning for the inner-loop in Variational Data Assimilation using Machine Learning}

\begin{document}
\maketitle
\begin{abstract}
\blfootnote{\sffamily Corresponding author: \email{victor.trappler@gmail.com}}
Data Assimilation is the process in which we improve the representation of the state of a physical system by combining information coming from a numerical model, real-world observations, and some prior modelling. It is widely used to model and to improve forecast systems in Earth science fields such as meteorology, oceanography and environmental sciences. One key aspect of Data assimilation is the analysis step, where the output of the numerical model is adjusted in order to account for the observational data. 
In Variational Data Assimilation and under Gaussian assumptions, the analysis step comes down to solving a high-dimensional non-linear least-square problem. In practice, this minimization involves successive inversions of large, and possibly ill-conditioned matrices constructed using linearizations of the forward model. In order to improve the convergence rate of these methods, and thus reduce the computational burden, preconditioning techniques are often used to get better-conditioned matrices, but require either the sparsity pattern of the matrix to inverse, or some spectral information.
We propose to use Deep Neural Networks in order to construct a preconditioner. This surrogate is trained using some properties of the singular value decomposition, and is based on a dataset which can be constructed online to reduce the storage requirements.
\end{abstract}
\begin{keywords}
  Variational Data Assimilation, Neural Networks, Preconditioning
\end{keywords}
\section*{Introduction}
Numerical models are ubiquitous nowadays as they are used to better understand and predict complex physical phenomena. In order to improve the accuracy and the predictability of those modelled systems, real-world data are assimilated into the predictions to provide a better representation of the true underlying state of the systems studied.
In Data Assimilation, this process is called the analysis step, where we combine different sources of information: the forecast coming from the previous time window, the available direct or indirect observations of various physical quantities within this time window, and some expert knowledge on the modelled processes, such as conservation and balance laws. 
Due to the time critical nature of those forecasts, the sheer size of the data involved, and the large computational power required to run numerical models, Data Assimilation methods have to be efficient since every improvement in those methods can lead to the use of more precise or more complex models, for a constant time budget.

In Variational Data Assimilation, the analysis is performed by minimizing a well-chosen objective function. This optimization can be very expensive since it happens in a high-dimensional space. Nonetheless, it can be tackled with gradient-based optimization, which boils down to successive high-dimensional linear system to solve. The speed of convergence of those methods depends on the condition number of the matrices involved, that is why several studies have been conducted on the condition number of various Data Assimilation problem, such as in \cite{haben_conditioning_2011,gurol_b_2014,tabeart_new_2021}.

Machine-Learning, on the other hand, has been increasingly applied on various aspects of Data Assimilation, as reviewed in \cite{cheng_machine_2023}. Some works focus on the  Data Assimilation process, as in \cite{boudier_dan_2020} where the authors propose a formalism of Data Assimilation, and apply recurrent Neural networks to perform the analysis and prediction steps. Same goes for \cite{arcucci_deep_2021}. In \cite{peyron_latent_2021}, an auto-encoder architecture is proposed in order to reduce the dimension of the state vector, and perform the assimilation in a lower dimensional latent space. Learning the underlying dynamical system is also of big interest. In~\cite{gottwald_combining_2021}, the authors propose to use Data Assimilation to learn the time-propagator of a dynamical system,  while in~\cite{dubois_data-driven_2020}, the whole dynamics of a Lorenz system is learned.

In this work, we propose to use Deep Neural Networks (DNN) to construct a preconditioner, not necessarily sparse, in order to improve the convergence of the Conjugate Gradient algorithm in a Variational Data Assimilation system. Using ML in Linear Algebra problems has recently found some traction in some related works:
in~\cite{ackmann_machine-learned_2021}, the authors build a preconditioner for an implicit solver, or in \cite{sappl_deep_2019,tang_graph_2022}, where a preconditioner for conjugate gradient is built using a convolutional neural network for the former, and a Graph Neural Network in the latter. In \cite{luna_accelerating_2021}, the authors proposes to use Neural Networks to improve the first guess in the GMRES method. Finally, in~\cite{hausner_neural_2023}, the authors manage to learn a sparse factorization of a matrix using Graph Neural Networks and the Frobenius norm, in order to precondition the Conjugate Gradient.

We will first review the classical method to obtain the inner/outer loop paradigm for optimization in order to introduce preconditioning, and then show how preconditioners can improve the convergence rate of CG, and how those can be constructed in a efficient way. 

\section{Variational Data Assimilation}
In what follows, we will first introduce the common notations used throughout this work and how the Variational Data assimilation process can be formulated as sequence of large-scale linear systems to solve.

\subsection{Data Assimilation as an optimization problem}
We assume that the physical system studied can be represented as a $n$-dimensional state vector $x \in \mathbb{X} \subseteq \mathbb{R}^n$. This state vector might represent different prognostic variables discretized on a mesh.
Let us consider a forward model $\mathcal{M}$ which maps the state-space onto itself. This operator usually represents the propagation in time of the state vector.
\begin{equation}
    \begin{array}{rcl}
    \mathcal{M}: \mathbb{X} \subseteq \mathbb{R}^n & \longrightarrow & \mathbb{X}  \\
        x & \longmapsto & \mathcal{M}(x) 
    \end{array}
\end{equation}
The output of the forward model (ie a state vector at a later time) often cannot be compared directly to the observations $y$. Indeed the observations may come from different sources, and are sparse and noisy quantities derived from the state. An observation operator $\mathcal{H}$ is then required to map the state vector to the observation space:
\begin{equation}
    \begin{array}{rcl}
    \mathcal{H}: \mathbb{X} \subseteq \mathbb{R}^n & \longrightarrow & \mathbb{Y} \subseteq \mathbb{R}^p  \\
        x & \longmapsto & \mathcal{H}(x)
    \end{array}
\end{equation}

In Data Assimilation, variational methods refer to approaches based on the optimization of an objective function, which measures the misfit between the model prediction and the observations, with a regularization that models the prior knowledge as a background term $x^b$ and $B$:
\begin{equation}
  J(x) =  \frac12\| \mathcal{G}(x) - y \|^2_{R^{-1}} + \frac12\| x - x^b\|_{B^{-1}}^2
  \label{eq:def_J}
\end{equation}
where the \emph{Generalized forward model} is
\begin{align}
  \mathcal{G}(x) &= (\mathcal{H} \circ \mathcal{M})(x)    \label{eq:G_def}
  \end{align}
  and the vector norms are defined for  $v \in \mathbb{R}^n$ and $\Sigma \in \mathbb{R}^{n\times n}$ positive definite as $\| v \|^2_{\Sigma} = v^T \Sigma v$.

From a probabilistic point of view, we can get to the same formulation by making the following Gaussian assumptions:
\begin{align}
y \mid x &\sim \mathcal{N}(\mathcal{G}(x), R) \\
x &\sim \mathcal{N}(x^b, B)
\end{align}
which leads to the expression of the objective function of Eq.~\cref{eq:def_J} as the negative log posterior probability of $x$ given $y$.




 \subsection{Incremental 4D-Var}
 In some large-scale systems, the Tangent Linear Model (ie the linearization of the model operator) and its adjoint may be available at the cost of proper derivation and maintenance, and at a computational cost roughly equivalent to the forward model. This means that we can consider gradient-based optimization methods in order to solve the analysis step.

Starting from a given state $x$, adding a small perturbation $\delta x$ gives
\begin{align}
  J(x + \delta x) &= \frac12 \| \mathcal{G}(x + \delta x) - y \|^2_{R^{-1}} + \frac12\| x + \delta x - x^b\|_{B^{-1}}^2
\end{align}
and linearizing $\mathcal{G}$ around $x$ gives the incremental version of the cost function
\begin{align}
  J_{\mathrm{inc}}(x, \delta x) &=  \frac12\|\mathcal{G}(x) + \G_x \delta x - y\|^2_{R^{-1}} + \frac12\|\delta x + x - x^b\|_{B^{-1}}^2 \\
   &= \frac12\| \G_x \delta x - d \|^2_{R^{-1}} + \frac12\|\delta x + x - x^b\|^2_{B^{-1}}
\end{align}
where $d = \mathcal{G}(x) - y$ are the departures from the observations and $\G_x = {\mathbf{H}}_{\mathcal{M}_x}\mathbf{M}_x$ is the Jacobian matrix of $\mathcal{G}$ evaluated at $x$.
Minimizing the incremental cost function with respect to $\delta x$ is a quadratic minimization problem, and the optimal increment $\delta x$ verifies 
\begin{equation}
\label{eq:inner_loop_linear_system}
  \underbrace{(\G_x^T R^{-1} \G_x + B^{-1})}_{\A_x} \delta x = \underbrace{-\G_x^T R^{-1} d - B^{-1}(x - x^{b})}_{b_x}
\end{equation}
and thus requires the resolution of a linear system of dimension $n$ using iterative methods, since the explicit inversion of such a matrix is unfeasible in practice.
A similar derivation can be achieved by applying Gauss-Newton Algorithm (see for instance \cite{gratton_approximate_2007}), which solves the original problem via successive approximations of the Hessian matrix of the non-linear optimization problem by the matrix $\A_x$.

One can also see the incremental formulation as a Bayesian Inverse Linear problem, where we are looking for the posterior mode (or posterior mean equivalently in this case) of $\delta x\mid d$
\begin{align}
    d \mid \delta x &\sim \mathcal{N}(\G_x \delta x, R) \\
    \delta x &\sim \mathcal{N}(x-x^b, B)
\end{align}
and the posterior mean is given by solving Eq.~\cref{eq:inner_loop_linear_system} and the posterior covariance matrix is
\begin{equation}
    \Gamma_{\text{post}} = \left(\G_x^T R^{-1} \G_x + B^{-1}\right)^{-1} = \A_x^{-1}
\end{equation}
Optimal approximations of this posterior are studied in \cite{benner_low-rank_2018,spantini_optimal_2015}.
\subsection{Nested loops}
Once the optimal increment $\delta x$ has been computed, the new point of linearization is chosen as $x + \delta x$, and a new approximation can be constructed. This can be repeated until convergence, or until a specified number of linearizations has been reached.

  \begin{figure}
     \centering
     \includegraphics[width=.4\textwidth]{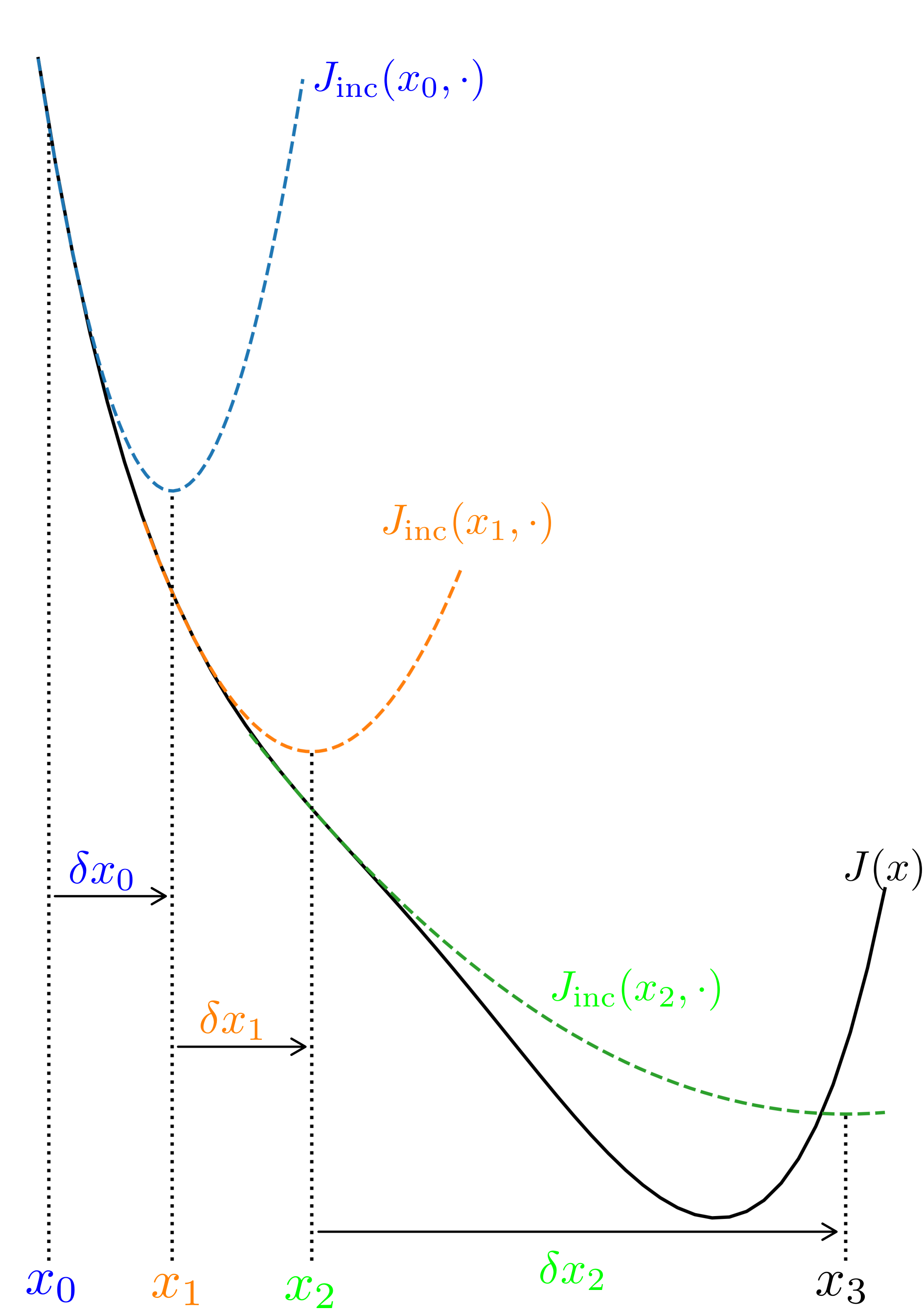}
     \caption{Illustration of minimization using successive quadratic approximations}
     \label{fig:succ_lin}
 \end{figure}
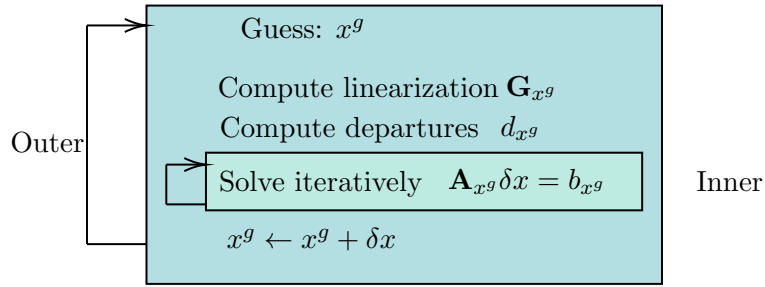
\begin{figure}
\begin{center}
      \scalebox{1.}{\tikzset{every picture/.style={line width=0.75pt}} 

\begin{tikzpicture}[x=0.75pt,y=0.75pt,yscale=-1,xscale=1]

\draw[fill=diagram!30]   (190,40) -- (450,40) -- (450,180) -- (190,180) -- cycle ;

\draw[fill=lightorange!30]   (220,114) -- (440,114) -- (440,143) -- (220,143) -- cycle ;
\draw    (200,140) -- (220,140) ;
\draw    (200,120) -- (200,140) ;
\draw    (200,120) -- (218,120) ;
\draw [shift={(220,120)}, rotate = 180] [color={rgb, 255:red, 0; green, 0; blue, 0 }  ][line width=0.75]    (10.93,-3.29) .. controls (6.95,-1.4) and (3.31,-0.3) .. (0,0) .. controls (3.31,0.3) and (6.95,1.4) .. (10.93,3.29)   ;
\draw    (160,50) -- (188,50) ;
\draw [shift={(190,50)}, rotate = 180] [color={rgb, 255:red, 0; green, 0; blue, 0 }  ][line width=0.75]    (10.93,-3.29) .. controls (6.95,-1.4) and (3.31,-0.3) .. (0,0) .. controls (3.31,0.3) and (6.95,1.4) .. (10.93,3.29)   ;
\draw    (160,50) -- (160,160) ;
\draw    (160,160) -- (190,160) ;

\draw (294.5,82.25) node  [xscale=1.0,yscale=1.0] [align=left] {{ Compute linearization}};
\draw (384.23,81.5) node  [xscale=1.0,yscale=1.0]  {$\mathbf{G}_{x^{g}}$};
\draw (231,45) node [anchor=north west][inner sep=0.75pt]  [xscale=1.0,yscale=1.0] [align=left] {{ Guess:}};
\draw (284,45) node [anchor=north west][inner sep=0.75pt]  [xscale=1.0,yscale=1.0]  {$x^{g}$};
\draw (290,102.25) node  [xscale=1.0,yscale=1.0] [align=left] {{ Compute departures}};
\draw (378.75,101.5) node  [xscale=1.0,yscale=1.0]  {$d_{x^{g}}$};
\draw (275.5,129.5) node  [xscale=1.0,yscale=1.0] [align=left] {{ Solve iteratively}};
\draw (382,128.5) node  [xscale=1.0,yscale=1.0]  {$\mathbf{A}_{x^{g}} \delta x = b_{x^{g}}$};
\draw (229,150.4) node [anchor=north west][inner sep=0.75pt]  [xscale=1.0,yscale=1.0]  {$x^{g} \gets  x^{g}  + \delta x$};
\draw (115.,102.27) node [anchor=north west][inner sep=0.75pt]  [xscale=1.0,yscale=1.0] [align=left] {{ Outer}};
\draw (461,122) node [anchor=north west][inner sep=0.75pt]  [xscale=1.0,yscale=1.0] [align=left] {{ Inner}};
\end{tikzpicture}}
\end{center}
\caption{Inner and outer loop paradigm for optimization}
\label{fig:inner_outer_schema}
\end{figure}

 The whole minimization procedure can be organized in nested loops, as detailed in \cref{fig:inner_outer_schema} and~\cref{alg:4dvar}.
\begin{itemize}
    \item The \emph{Outer Loop}, which requires a run of the forward model $\mathcal{G}$ at a point $x$, and the evaluation of the Tangent Linear Model $\G_{x}$ in order to get a linearization. The linearization of the cost function, which implies the Tangent Linear Model, can be obtained by classical methods of automatic differentiation. The number of outer loops is critical when dealing with highly non-linear processes (\cite{bonavita_nonlinear_2018})
    
  \item  the \emph{Inner Loop}, where we solve the minimization problem using the TLM (ie successive quadratic approximations). Once this minimization has been performed, the point of evaluation for the Outer Loop is chosen.
\end{itemize}


\begin{algorithm}
\begin{algorithmic}
\State $n \gets 1$
\State $x_i \gets x_0$
\While{$i \leq n_{\text{outer}}$} \Comment{\textbf{Outer Loop}}
\LComment{Direct model and linearization at $x_i$}
    \State Evaluate $\mathcal{G}(x_i)$, $J(x_i)$, $\G_{x_i}$
    \State $b_{x_i} \gets -\G_{x_i}^TR^{-1}\left(\mathcal{G}(x_i) - y\right)$
    \State $\A_{x_i} \gets (\G_{x_i}^T R^{-1} \G_{x_i} + B^{-1})$ 
    \LComment{The linear system to solve is $\A_{x_i} \delta x_i = b_{x_i}$}
    \State $j \gets 0$, $\delta x^{(j)} \gets 0$
    \While{$j \leq n_{\text{inner}}$ or $\|r_j\|_2 < \epsilon$} \Comment{\textbf{Inner Loop}}
        \State $\delta x^{(j+1)} \gets \text{ConjugateGradient}(\A_{x_i}, \delta x^{(j)})$
        \State $r_j \gets \A_{x_i}\delta x^{(j+1)} - b_{x_i}$
        \State $j\gets j+1$
    \EndWhile
    \State $\delta x_i \gets \delta x^{(j)}$
    \State $x_{i+1} \gets x_i + \delta x_i$
    \State $i \gets i+1$
\EndWhile
\end{algorithmic}
\caption{Pseudocode of the minimization procedure in 4DVar}
\label{alg:4dvar}
\end{algorithm}


\subsection{Conjugate Gradient}
In the inner loop, the matrix $\A_x$ cannot be constructed explicitely, let alone be inverted via direct methods. We can use Krylov-subspace based methods to approximately solve the linear system, such as GMRES, or Conjugate Gradient which only require matrix-vectors products.
Since the matrix to inverse is symmetric positive definite (spd), we use the Conjugate Gradient algorithm to solve the linear system (see \cite{freitag_numerical_2020, diouane_efficient_2024} for specifics of CG in Data Assimilation) and the error $e_k = \delta x_k - \delta x^*$ between the computed increment at the $k$th step and the true value $\delta x^* = \A_x^{-1} b_x$ can be bounded, giving a rough rate of convergence
\begin{align}
\label{eq:cond_number}
    \|e_k\| &\leq 2\left(\frac{\sqrt{\kappa(\A_x)} - 1}{\sqrt{\kappa(\A_x)} + 1}\right)^k \|e_0\| 
\end{align}
where
 $\kappa(\A_x) = \|\A_x^{-1}\|_2\cdot \|\A_x\|_2 \geq 1 = \kappa(I_n)$ is the condition number of the matrix $\A_x$. As this matrix is symmetric positive definite, this condition number can be written as the ratio between the largest and smallest eigenvalues:
 \begin{equation}
     \kappa(\A_x) = \frac{\lambda_1(x)}{\lambda_n(x)}
 \end{equation}
 where the spectrum of $\A_x$: $\mathrm{sp}(\A_x)=(\lambda_1(x),\dots,\lambda_n(x))$ is sorted in descending order.
 
 It is clear from Eq.~\cref{eq:cond_number} that a condition number close to $1$ leads to a better convergence rate of the CG algorithm. Since the matrix $\A_x$ is fully determined by the problem, its condition number is not directly adjustable. We can however use a preconditioner in order to improve the condition number of the problem, and thus improve the convergence rate for this iterative method.

 \subsection{Preconditioning the Inner Loop}
Instead of directly solving the linear system  $\A_x \delta x = b$ using iterative methods, one can look for a system which possesses the same solution, ie $\A_x^{-1}b$, but for which the CG method converges faster. One approach is to left multiply the two sides of the equation by an invertible matrix of size $n\times n$, say $L^T$ giving the linear system $(L^T\A_x) \delta x = (L^Tb)$.

In order to conserve the symmetric property of the matrix to inverse and use CG, we can rewrite the linear system as
   \begin{align}
    \underbrace{(L^T \A_x L)}_{\tilde{\A}}\underbrace{(L^{-1} \delta x)}_{\tilde{x}} &= L^T b
\end{align}
If $\tilde{x} \in \mathbb{R}^n$ verifies the linear equation $\tilde{\A} \tilde{x} = L^Tb$, the solution of the original linear system can be retrieved by $\delta x = L \tilde{x}$. 
The new linear system can also be preconditioned if needed, but we focus here on "first-level" preconditioning.

The matrix $P=LL^T$ is called a preconditioner, while $L$ is sometimes called a split preconditioner, and
$P\A_x$ and $L^T\A_x L$ share the same spectrum.
Trivial examples of preconditioners include $P=I_n$ and $P=\A_x^{-1}$, but for the former the problem to solve is left unchanged, while for the latter the solution is found trivially, at the cost of computing directly the inverse of the matrix.
The choice of a preconditioner is largely problem dependent, but some desirable properties can be listed:
\begin{itemize}
    \item $P$ should be symmetric and non-singular
    \item $P$ should be cheap to apply as a linear operator
    \item $P$ should improve the condition number of $\A_x$ in order to improve the convergence of iterative methods
\end{itemize}
In data assimilation, given the definition of $\A_x$ in Eq.~\cref{eq:inner_loop_linear_system} , particular choices of $L$ can be useful to simplify the problem. Indeed, preconditioning the matrix $\A$ using $L=B^{-1/2}$ gives 
\begin{equation}
    \tilde{\A} = B^{-T/2}\G_x^TR^{-1}\G_xB^{-1/2} + I_n
\end{equation}
In this case, all the eigenvalues of $\tilde{\A}$ are larger than $1$, so its condition number is smaller than its largest eigenvalue (see \cite{gurol_b_2014}).

In many cases, one may look for a solution of the linear system in a smaller subspace generated by the columns of $L$. This method is often named in the literature Control Variable Transform, and thus $L$ is not a square matrix. However the two problems are not necessarily equivalent, and~\cite{menetrier_overlooked_2015} studies further the conditions for equivalence. In the case of sparse matrices, a preconditioner can be found by looking for a product $P\A_x$ which approximates the identity matrix. That is the principle of Sparse Approximate Inverse~(see \cite{grote_parallel_1997}), where the preconditioner is found by minimizing $\|I_n - P\A_x\|$ for $P$ with a prescribed sparsity pattern.

Since the convergence properties of the CG method is dependent on the distribution of the eigenvalues of the matrix $\A_x$, we will focus on preconditioners constructed using  its spectral properties.


\subsection{Spectral preconditioners}
We will now drop the subscript $x$ for notation sake, but all those quantities depend implicitely on the point of linearization $x$.
The main idea behind spectral preconditioners is to act directly on the eigenvalues of $\A$, by constructing a matrix which will decrease the $r$ largest eigenvalues of $\A$ to some smaller values, thus decreasing the ratio defining the condition number in Eq.~\cref{eq:cond_number}.
The spectral preconditioners introduced here are studied more generally as Limited Memory Preconditioners in \cite{tshimanga_limited-memory_2008}.

Since $\A$ is symmetric positive definite, eigendecomposition and singular value decomposition are equivalent.
Let $\A=U\Lambda U^T$ be the Singular Value Decomposition (SVD) of $\A$ with $U = (u_1 \mid u_2 \mid \dots \mid u_n) \in \mathbb{R}^{n \times n}$ an orthonormal matrix, and $\Lambda = \mathop{\mathrm{diag}}(\lambda_1,\dots,\lambda_n)$ 
where the $\lambda_i$ are all strictly positive and sorted in descending order.

Truncating the SVD on its $r$ first components gives the low-rank approximation of $\A$:
\begin{align}
\label{eq:def_lr}
   \A_r = U_r\Lambda_r U_r^T
\end{align}
where $U_r = \left(u_1 \mid \dots \mid  u_r\right) \in \mathbb{R}^{n \times r}$, and $\Lambda_r = \mathrm{diag}(\lambda_1, \dots, \lambda_r)$.

Eckart–Young–Mirsky theorem provides another characterization of the low-rank approximation, in terms of an optimization problem, which will be used in \cref{ssec:loss_def}:
\begin{equation}
\label{eq:eym_thm}
    \min_{\tilde{\A}; \text{rk}(\tilde{\A}) = r} \| \A - \tilde{\A} \|^2_{\text{F}} =\| \A - \A_r \|^2_{\text{F}} = \sum_{i=r+1}^n \lambda^2_{i}
\end{equation}
where $\|\cdot \|_{\text{F}}$ is the Frobenius matrix norm defined for a matrix $D$ as
\begin{align}
    \|D\|_{\text{F}}^2 &= \mathrm{tr}\left(DD^T\right) = \sum_{i,j} d_{ij}^2 
\end{align}

The different terms of the decomposition of Eq.~\cref{eq:def_lr} can be used to construct a symmetric matrix $\precmat{\alpha}$ that will act as a preconditioner by treating differently the leading eigenvalues and the remaining ones:
\begin{equation}
\label{eq:def_Palpha_exact}
    \precmat{\alpha} = \beta I_n + U_r(\mu\Lambda_r^{\alpha} - \beta I_r)U_r^T
\end{equation}
with
\begin{itemize}
\item $\alpha$ the exponent of the eigenvalues to consider
    \item $\mu > 0$ the value which will affect the $r$ leading eigenvalues
    \item $\beta > 0$ the value which will multiply the $n-r$ other eigenvalues
\end{itemize}
This type of scaled preconditioners have been studied more thoroughsly from a theoretical point of view in  \cite{diouane_efficient_2024}.
We can better understand the effect of this matrix on an arbitrary vector $x \in \mathbb{R}^n$ by decomposing it into an element  $x_r \in \mathrm{range}(U_r)$, the span of the first $r$ eigenvalues, and an element $x_\perp$ in its null-space. 
There exists then $w \in\mathbb{R}^r$ such that $x = x_r + x_\perp = U_rw + x_\perp$. Applying $P_\alpha$ gives
\begin{align}
    \precmat{\alpha}x &= (\beta I_n + U_r(\mu\Lambda_r^{\alpha} - \beta I_r)U_r^T)(U_rw + x_\perp) \\ \nonumber
    &= U_r\left(\mu \Lambda_r^{\alpha}\right)w + \beta x_\perp
\end{align}
thus the components in $\mathrm{range}(U_r)$ are multiplied by the diagonal matrix $\mu \Lambda_r^\alpha$, while the components in the null-space are multiplied by $\beta$. 

\begin{figure}
    \centering
    \includegraphics[scale=0.5
]{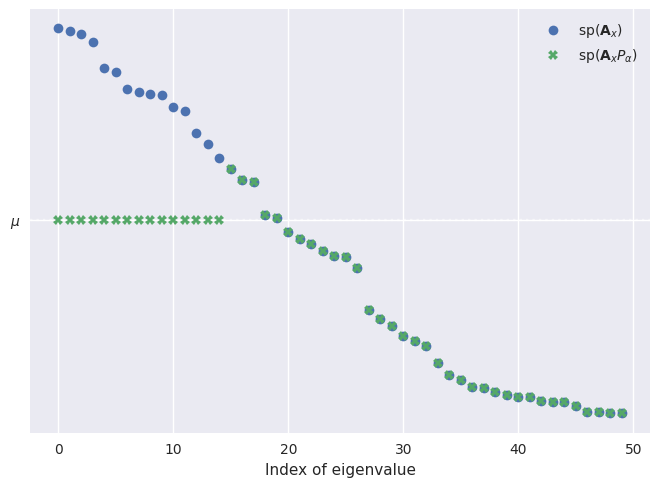}
    \caption{Illustration of the spectrum of an example spd matrix $\A_x$, and the preconditioned matrix using $P_{\alpha}$ for $\alpha=-1$ }
    \label{fig:spectrum}
\end{figure}

By construction, $\precmat{\alpha}$ is a spd matrix with spectrum 
\begin{equation}
    \mathrm{spectrum}(\precmat{\alpha}) = \{\mu\lambda_1^\alpha, \dots , \mu\lambda_r^\alpha, \beta\dots, \beta\}
\end{equation}
and $\precmat{\alpha/2}$ is a matrix square root of $\precmat{\alpha}$.
Since $\A$ and $\precmat{\alpha}$ share the same eigenvectors, the spectrum of the product is
\begin{equation}
         \mathrm{spectrum}(\precmat{\alpha/2}^T\A \precmat{\alpha/2}) =  \mathrm{spectrum}(\A \precmat{\alpha}) = \{\mu\lambda_1^{\alpha+1} \dots , \mu\lambda_r^{\alpha+1}, \beta\lambda_{r+1}\dots, \beta\lambda_{n}\}
\end{equation}
This spectrum highlights how to set $\alpha$ and $\beta$ to construct a preconditioner: choosing $\alpha = -1$ and $\beta=1$, as in \cref{fig:spectrum}, groups the $r$ leading eigenvalues of the matrix product to $\mu$, so by choosing $\mu$ inbetween the smallest eigenvalue $\lambda_n$ and $\lambda_r$, the condition number of the preconditioned matrix $\A \precmat{\alpha}$ is less than $\frac{\lambda_{r}}{\lambda_n}$.

Such a preconditioner can be used to cluster the $r$ leading eigenvalues at $\mu$, and thus improve the convergence rate in the Conjugate Gradient algorithm. However, a precise computation of the SVD might be challenging in practice: methods such as the Lanczos iterations require the evaluations of many matrix-vector products (usually more than $r$). Recently, randomized methods have been proposed for these kind of computations in data assimilation, see~\cite{dauzickaite_randomised_2021}.
Those procedures are dependent on the matrix $\A=\A_x$ at the point of linearization $x$, so even if some eigen-information can be reused when the linearization point does not change much, as done in \cite{tshimanga_limited-memory_2008}, most computations are discarded at the start of a new assimilation window.

Instead, we propose to use Deep Neural Networks in order to map the state of linearization $x$ to an approximate low-rank decomposition of $\A_x$ which can be used as a preconditioner.

\section{Deep Neural Network to construct state-dependent preconditioners}
\subsection{Architecture of the Deep Neural Network}
In order to construct a preconditioner based on Eq.~\cref{eq:def_Palpha_exact}, two elements are needed: we need to approximate $U_r$ by a matrix of size $n\times r$, whose columns are orthonormal, and a vector of size $r$, with positive elements to approximate $\Lambda_r$. 

We propose to use a Deep Neural Network (parameterized by $\theta$, a vector containing all the weights and biases of this DNN), say $f_\theta$, in order to compute those to produce tensors of appropriate dimensions. Given $x\in \mathbb{R}^n$, this Neural Network outputs both a set of $r$ non-orthonormal vectors $\tilde{U}_\theta(x)\in \mathbb{R}^{n \times r}$, and a vector $\tilde{\Lambda}_\theta(x)\in\mathbb{R}^{r}$, which are to be postprocessed in order to verify the aforementioned properties. This allows for a flexible choice of the architecture of the DNN, which can then be chosen in a problem specific manner (CNN for spatially distributed states for instance).

\begin{figure}
    \centering
    \includegraphics[width=0.6\textwidth]{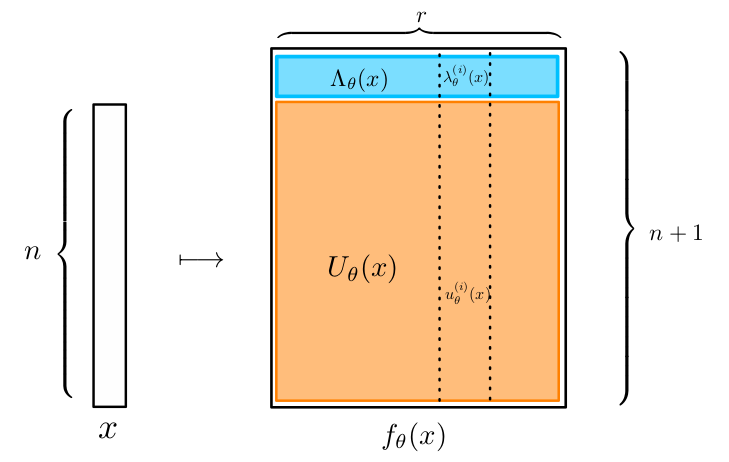}
    \caption{Schematic representation of the input/output signature of the Neural Network}
    \label{fig:in_out_nn}
\end{figure}

To ensure the orthonormal property of the vectors, we use the QR decomposition on $\tilde{U}_\theta(x)$, which is numerically stable compared to a classical Gram-Schmidt orthonormalization procedure, while the positivity of the approximate eigenvalues $\tilde{\Lambda}_\theta(x)$ is imposed using any function $\mathbb{R} \rightarrow \mathbb{R}^+$ elementwise. In this work, we will use a scaled sigmoid function: $S_M: x \longmapsto \frac{M}{1+e^{-x}}$, where $M$ can be chosen as a rough upper bound on the singular values of $\A_x$. This choice allows for bounding the resulting eigenvalues into an acceptable range, which helps avoid numerical issues during training.
This mapping is summarized Eq.~\cref{eq:NN_flow}, and \cref{fig:in_out_nn}.
\begin{equation}
    \begin{array}{rclclcl}
         \mathbb{R}^n& \longrightarrow &\mathbb{R}^{(n + 1)r}& \longrightarrow&\mathbb{R}^{n \times r} \times \mathbb{R}^r & \longrightarrow&\mathbb{R}^{n \times r} \times \mathbb{R}_+^r\\
         x& \overset{\text{DNN}}{\longmapsto} & f_{\theta}(x)&\overset{\text{split}}{\longmapsto} & (\tilde{U}_{\theta}(x), \tilde{\Lambda}_{\theta}(x))&\overset{\texttt{qr},S_M}{\longmapsto}& \left(\texttt{qr}(\tilde{U}_{\theta}(x)), S_M\left(\tilde{\Lambda}_{\theta}(x)\right)\right)  \\
         & & & & & & =\left(U_\theta(x), \Lambda_\theta(x)\right)
    \end{array}
    \label{eq:NN_flow}
\end{equation}

Given the output of $f_\theta$ and a postprocessing using QR and $S_M$, the Neural Network-based low-rank reconstruction of rank $r$ is
\begin{align}
\label{eq:def_ml_lr_approx} 
\A_\theta(x) &= U_{\theta}(x) \Lambda_{\theta}(x) U_{\theta}(x)^T \\
&= \sum_{i=1}^r \lambda_{\theta}^{(i)}(x) u^{(i)}_{\theta}(x)\label{eq:def_ml_rk1_approx} 
\left(u^{(i)}_{\theta}(x)\right)^T
\end{align}
with $U_{\theta}(x) = \left(u^{(1)}_\theta(x)\mid \dots \mid u^{(r)}_\theta(x)\right)$. Using the decomposition, a split preconditioner can be defined as using Eq.~\cref{eq:def_Palpha_exact} for $\beta=1$ and $\alpha = -1/2$:
\begin{equation}
\label{eq:def_split_prec_ml}
    L_\theta(x) =  I_n + U_{\theta}(x)\left(\mu \Lambda_{\theta}(x)^{-1/2} - I_r\right)U_{\theta}(x)^T \quad \text{ with } \mu \geq 1
\end{equation}

In theory, if the DNN provides the optimal low-rank approximation of $\A_x$, choosing $\mu=1$ would allow to group all the $r$ first eigenvalues to $1$, thus reducing the condition number of the matrix. In practice, the DNN only produces an approximation of the eigenvectors and of the eigenvalues, meaning that there is a risk to worsen the condition number. Experiments have shown that choosing $\mu \in [\min_i \lambda_{\theta}^{(i)}(x), \max_i \lambda_{\theta}^{(i)}(x)]$ helps to account for the approximation error due to the DNN. This is further discussed in \cref{sec:sw_application}.
\subsection{Loss function definition using Frobenius norm approximation}
\label{ssec:loss_def}
Neural networks are parameterized by $\theta \in\mathbb{R}^{\mathfrak{N}}$, which combines all the weights and biases of the individual neurons of $f_\theta$. To set this parameter, one need to define an appropriate metric which is then optimized. 
Given the Eckart–Young–Mirsky theorem Eq.~\cref{eq:eym_thm}, which defines the SVD in terms of an optimization problem and the reconstruction defined in Eq.~\cref{eq:def_ml_lr_approx}, we define the loss for a single state of linearization $x_i$ as
\begin{equation}
\label{eq:loss_per_point}
    \mathcal{L}_{\text{explicit}}(\theta; x_i)=\|\A_\theta(x_i) -\A_{x_i}\|^2_{\text{F}}
\end{equation}
where this term would be minimal if $\A_{\theta}(x_i)$ is the low-rank approximation of $\A_{x_i}$.

This loss requires the evaluation of the norm of the difference of two $n\times n$ non-sparse matrices, which brings several challenges. Constructing the matrix $\A_{x_i}$ is computationally expensive, since 
in most differentiated computer codes, this matrix is only accessible as an operator. In Data Assimilation especially, given the definition of $\A_x$ in Eq.~\cref{eq:inner_loop_linear_system}, computing $\delta x \mapsto \A_x \delta x$ requires the applications of two linear (with respect to the second argument) operators:
The Tangent Linear operator
\begin{equation}
    \mathrm{TL} : (x_i, \delta x) \longmapsto \G_{x_i}\cdot \delta x 
    \end{equation}
    and the adjoint operator
    \begin{equation}
    \mathrm{Adj}: (x_i, y) \longmapsto \G_{x_i}^T \cdot y
\end{equation}
From a computational point of view, applying one of those operators is within the same order of magnitude of complexity as the forward model $\mathcal{G}$.
Obviously, in order to construct the full Jacobian matrix $\A_x$, one could apply the linear operator to each vector $e_i$ of the canonical basis since $A_{x} = \left(A_x e_1 \mid \dots, ,\mid A_x e_n\right)$, but this is impractical since it requires $n$ evaluations, on top of the large memory requirements needed to store the matrix $\A_x$ for a single linearization point.

Same goes for the matrix $\A_\theta(x)$: constructing the full matrix is hard from a storage point of view, even though using it as a linear operator is cheaper since it requires only $r$ dot products of $n$-dimensional vectors as seen from Eq.~\cref{eq:def_ml_rk1_approx}, 

Since we are only interested in the Frobenius norm of the difference of the operators, we can instead directly estimate it using statistical estimators.
Let $D$ be a real matrix of size $n\times n$. Its squared norm $\|D \|^2$ can be rewritten as the expectation of a vector norm using the linearity of the trace and expectation operator:
\begin{align}
    \Ex_{\xi}\left[\|D \xi\|^2\right] = \Ex_{\xi}\left[\mathrm{tr}\left(D\xi\xi^TD^T\right)\right]=\mathrm{tr}\left(\Ex_{\xi}\left[\xi\xi^T\right]D^TD\right) = \|D \|_\text{F}^2
\end{align}
where $\xi\sim \mathcal{N}(0, I_n)$.
Given a matrix $Z\in\mathbb{R}^{n \times k}$ whose $k$ columns $z^{(j)}$ are sampled from a standard Gaussian distribution, we can use a Monte-Carlo estimator of the expectation:
\begin{equation}
\label{eq:estim_frob_norm}
    \frac{1}{k}\|DZ\|^2_{\text{F}}=\frac{1}{k}\sum_{j=1}^k \|Dz^{(j)} \|^2 \quad \text{estimator of} \quad \|D\|^2_{\text{F}}
\end{equation}
Other estimators of this norm using random samples are studied in~\cite{gudmundsson_small-sample_1995,gratton_improved_2018}, while in \cite{indyk_learning-based_2019}, the authors use ML to construct the matrix to evaluate. 

Using Eq.~\cref{eq:estim_frob_norm}, for a state-vector $x_i$ in the training dataset and $z_i^{(1)},\dots z_i^{(k)}$ i.i.d. samples of a standard Gaussian random variable, an estimate of the matrix norm of Eq.~\cref{eq:loss_per_point} is
\begin{align}
\label{eq:def_loss}
    \mathcal{L}_{\text{explicit}}(\theta; x_i) \approx \mathcal{L}(\theta; x_i) = \frac{1}{k}\sum_{j=1}^k \|{\A}_\theta(x_i)z_i^{(j)} - \A_{x_i}z_i^{(j)}\|^2
\end{align}
where ${\A}_\theta(x_i)$ is defined as in Eq.~\cref{eq:def_ml_lr_approx}.
We can also use the same estimator in order to estimate the norm of $\A_x$ as $\frac{1}{k}\sum_{j=1}^k \|\A_x z^{(j)}\|^2$, which is an estimate of the sum of all its eigenvalues squared. This can be used in order to normalize the loss in Eq.~\cref{eq:def_loss}, and can be interpreted as the fraction of \emph{unexplained} variance, by analogy with classical Principal Components Analysis:
\begin{equation}
    \mathcal{L}_{\text{relative}}(\theta, x_i) = \frac{\mathcal{L}(\theta;x_i)}{\frac{1}{k}\sum_{j=1}^k \|\A_{x_i} z^{(j)}_i\|^2}
\end{equation}

\subsection{Construction and storage of the training dataset}
In order to train the Neural Network, the construction of a dataset is needed in order to optimize the loss function defined in Eq.~\cref{eq:def_loss}. Each element (indexed by $i$) in this dataset consists of three elements: a state  $x_i$ which is used for the linearization, a random matrix $Z_i=(z_i^{(1)} \mid \dots \mid z_i^{(k)})\in \mathbb{R}^{n\times k}$ whose components are iid and normally distributed, and finally the evaluation of this sample by the matrix of interest: $\A_{x_i}Z_i$. The training dataset is then
\begin{align}
\label{eq:dataset_training}
    \mathfrak{D}_{\text{training}} = \left\{\left(x_i, Z_i, \A_{x_i}Z_i)\right) \in \mathbb{R}^{n} \times \mathbb{R}^{n\times k} \times \mathbb{R}^{n\times k}\quad \text{s.t.} \quad 1 \leq i \leq N_{\text{training}}\right\}
\end{align}
However, we do not have to store all the training set in memory: $Z_i$ is independent of $x_i$, and can be sampled when needed, and $\A_{x_i}$ depends only on $x_i$. 

\begin{algorithm}
\begin{algorithmic}
\For{$1\leq i \leq n_\text{batch}$}
\State Sample and store $Z_i = (z_i^{(1)}\,  | \dots | \, z_i^{(k)}) \in \mathbb{R}^{n \times k}$ with $z_{i}^{(j)} \sim \mathcal{N}(0, I_n)$ iid for $1\leq j \leq k$
\State Compute and store $\A_{x_i}Z_i\in \mathbb{R}^{n \times k}$
\State $x_{i+1} \gets $ New state generated from $x_i$
\EndFor
\end{algorithmic}
\caption{Pseudocode for the generation of a batch for online training}
\label{alg:iterabledataset}
\end{algorithm}

The method to generate a batch of $n_\text{batch}$ samples is summarized~\cref{alg:iterabledataset}. In order to train a Deep Neural Network, the constructed batches should be representative enough of the whole state space. To get appropriate diversity in the states used to build the batch, we propose to generate the new state iteratively by advancing the current state using the numerical model $\mathcal{M}$ with a randomly generated lead time,  large enough so that the $x_i$ used for the batch are not too correlated, and by potentially adding a small random perturbation before propagation.

\section{Application to a Shallow Water Assimilation system}
\label{sec:sw_application}
\subsection{Shallow Water equations and Data Assimilation setting}
The Shallow Water equations describe the motion of large bodies of water, for which the horizontal scale is larger than the vertical scale which is the case for rivers, seas and oceans. They consist in PDEs obtained by vertically averaging the Navier-Stokes equations. In this application, the variables of interest are the deviation of sea surface height $\eta$ around a mean height $\eta_0$, the velocity $u$ in the $x$-direction, and $v$, the velocity in the $y$-direction.
\begin{equation}
\left\{
    \begin{array}{rl}
        \frac{\partial \eta}{\partial t} + \frac{\partial (\eta_0 + \eta)u}{\partial x} + \frac{\partial (\eta_0 + \eta)v}{\partial y}&=0  \\
        \frac{\partial u}{\partial t} - \xi v + \frac{\partial B}{\partial x}&= \nu \Delta u - c_b u  + \frac{\tau_x}{\rho_0\eta_0} \\
        \frac{\partial v}{\partial t} + \xi u + \frac{\partial B}{\partial y} &= \nu \Delta v - c_b v  
    \end{array}
    \right.
\end{equation}

Those equations are discretized using a Arakawa C-grid of $64 \times 64$ cells, on a square domain of size $L_x = L_y = 1800\si{\kilo \meter}$, meaning that the three prognostic variables are $\eta \in \mathbb{R}^{64 \times 64}$, $u \in \mathbb{R}^{63 \times 64}$ and $v \in \mathbb{R}^{64 \times 63}$.
Once flattened and concatenated, the state vector is then $x =(\eta, u, v) \in \mathbb{R}^{12160}$. Explicitely storing the Gauss-Newton matrix would require $4.7 \si{\giga\byte}$ (without exploiting the symmetry).

We consider the model $\mathcal{M}$ that simulates the evolution of the state vector with a lead time of $T$ corresponding to 2 days.
\begin{equation}
    \begin{array}{rcl}
      \mathbb{R}^n   & \longrightarrow&\mathbb{R}^n  \\
        \mathcal{M}: x_t & \longmapsto &  \mathcal{M}(x_t) = x_{t + T} 
    \end{array}
\end{equation}

The cost function is defined as in Eq.~\cref{eq:def_J}
\begin{equation}
    J(x) = \frac{1}{2} \|(\mathcal{H}\circ \mathcal{G})(x) -y \|_{R^{-1}}^2 + \frac{1}{2}\|x - x^b \|_{B^{-1}}^2
\end{equation}
where $\mathcal{H}(x) = \mathcal{H}((\eta, u, v)) = \eta$, $R = I_{64^2}$, meaning that only the free-surface height is observed.
The background state $x^b \in \mathbb{R}^{n}$ is computed as the average of states obtained during a previous simulation with a large lead time.

\subsection{Neural Network Architecture}
For this problem, the state vector represents three spatial variables, arranged on a regular grid. By padding the $u$ and the $v$ component, we can reshape the state vector as a tensor of shape $(64, 64, 3)$, ie like an image with 3 channels. Each of those components is scaled so that each channel has approximately unit variance.
Because of this image-like structure, we can use Neural Network architecture well-suited for such data, such as Convolutional Neural Networks (CNN) or U-Nets. We found that using a U-Net architecture, with transformers instead of CNN for the subsampling step has shown good results for this problem.

\subsection{Dataset and training}
The training dataset is constructed according to Eq.~\cref{eq:dataset_training}, where $N_{\text{training}} = 1000$ states of linearization have been sampled, and $k=100$ random vectors have been used for matrix-vector products.

\subsection{Numerical Results}

We trained a DNN whose architecture allows us to get $r_{\text{train}}=2000$ approximate singular vectors and values, sorted by descending singular value.
Based on this, we can compare the preconditioners obtained using a different numbers of retained vectors (denoted as "rank", even though $\precmat{\alpha}$ is full-rank) : $r=1000$,  $1900$ and $2000$. For each of those, different values of $\mu$ have been chosen: either it is set to a fixed value, or it is set to the smallest eigenvalue provided by the DNN.
For $r=1000$, $\min_i\lambda^{(i)}\approx 160$, for $r=1900$, $\min_i\lambda^{(i)}\approx 145$, and finally, for $r=2000$, $\min_i\lambda^{(i)}\approx 2$.
The matrices to inverse have their leading eigenvalues close to 20000, and show approximatively an exponential decay.

In order to compare numerical results, we started from a base state $x^{\mathrm{base}}$. We generated the "truth" by perturbating and advancing the base state using the numerical model a random number of time steps. Is is then used to generate observations using $\mathcal{H}$ and by adding an observation noise.
The state of linearization ($x_i$) is chosen in a similar way, by perturbating the base state.
\begin{algorithm}
\begin{algorithmic}
\State{$x^{\dagger} \gets x^{\mathrm{base}} + \text{small perturbation}$}
\State{$x_1 \gets x^{\mathrm{base}} + \text{small perturbation}$}
\For{$1\leq i \leq N$}
\LComment{Generate Observations}
\State Advance the truth $x^{\dagger}$ by a random number of time steps
\State Compute observations $y \gets \mathcal{H}(x^{\dagger}) + \epsilon$
\LComment{Form the Linear system}
\State Linearize the forward model at $x_i$ to form the GN linear operator $\A_{x_i}$
\State Compute $b_{x_i}$ using the departures $\mathcal{G}(x_i) - y$ and the background
\LComment{Preconditioned CG}
\State Use the DNN to compute $L_\theta(x_i)$
\State Use CG for the linear system  $\A_{x_i}, b_{x_i}$ with $L_{\theta}(x_i)$ as preconditioner, to get $\delta x_i$
\LComment{Modify the linearization step}
\State $x_i \gets x_i + \delta x_i$
\EndFor
\end{algorithmic}
\caption{Pseudocode for the numerical experiment}
\label{alg:num_experiment}
\end{algorithm}

 Since we are solving iteratively a system of the form $\A_x \delta x = b$ (the subscript $i$ is dropped for convenience), the quantity of interest chosen to track the convergence of the Conjugate Gradient method is often the $L_2$ norm of the residual $e_j = \A_x\delta x^{(j)} - b$.
However, the CG method does not guarantee a monotonic decrease of the Euclidian norm of the residuals $\|e_j\|_2$, nor its energy norm $\|e_j\|_A$, which can explain some oscillations in some visualizations.

\begin{figure}[h]
    \centering
    \includegraphics[width=0.8\textwidth]{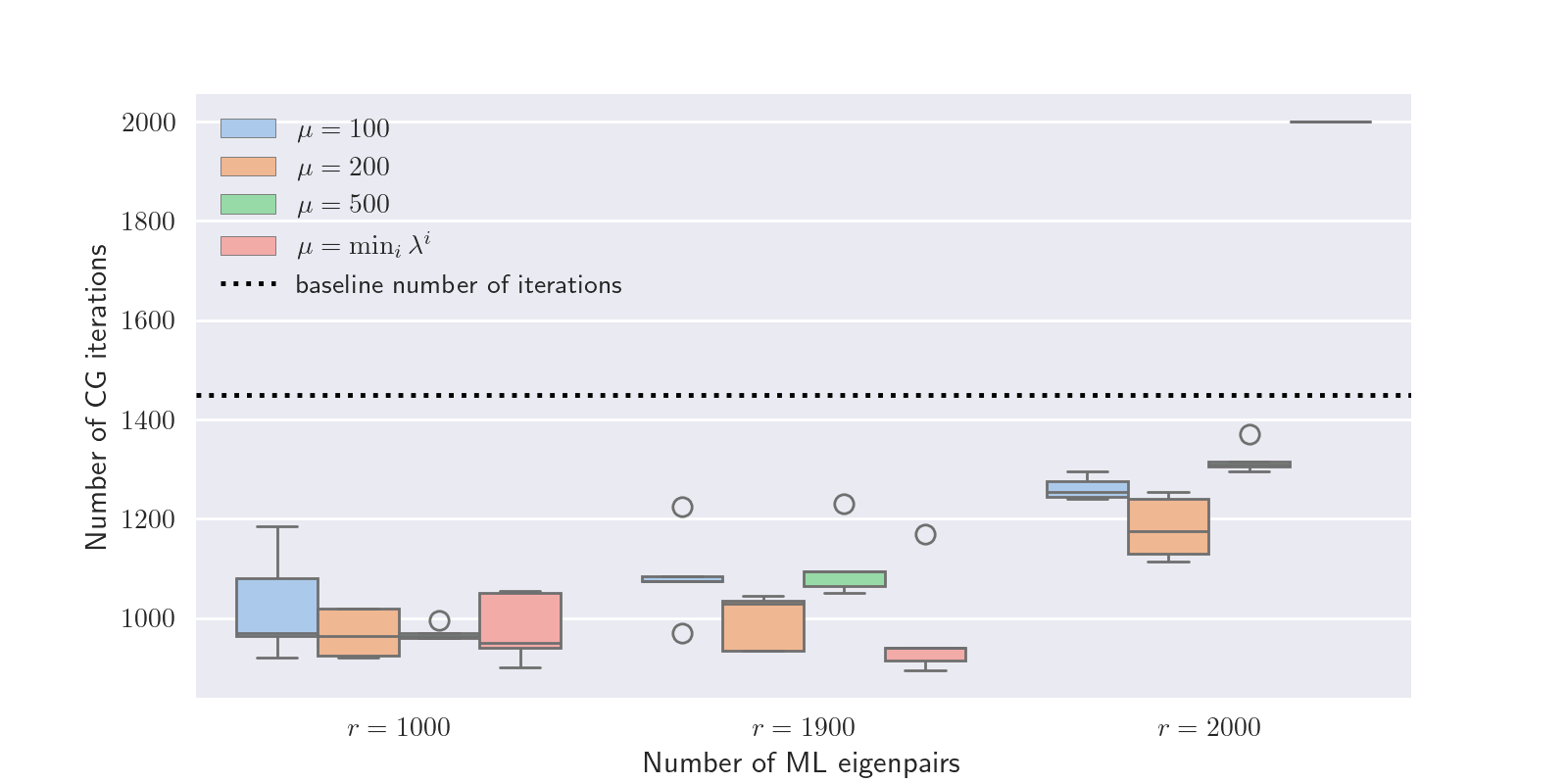}
    \caption{Number of iterations needed to reach the norm threshold. Dotted line indicates the number of iterations for the unprecondition problem.}
    \label{fig:boxplot_niters}
\end{figure}
\begin{figure}
    \centering
    \includegraphics[width=\textwidth]{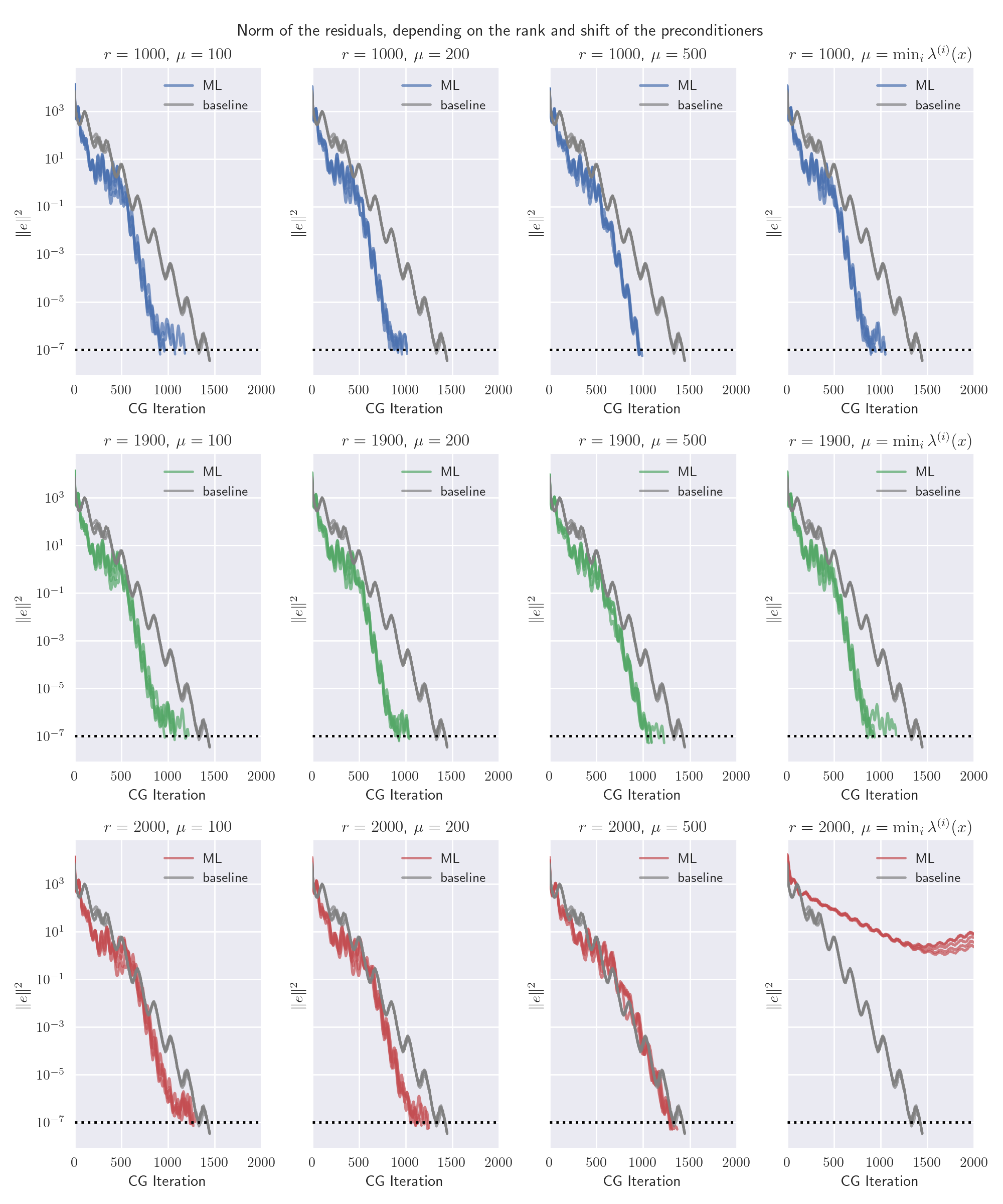}
    \caption{Norm of the residuals as a function of the iteration number during CG, depending on the number of retained vectors $r$ and shift parameter $\mu$ as defined in Eq.~\cref{eq:def_split_prec_ml}}
    \label{fig:norm_res}
\end{figure}

\Cref{fig:boxplot_niters} shows the number of iterations needed to reach the threshold of $10^{-7}$ for the euclidian norm of the residuals or when 2000 iterations of CG has been reached (whichever comes first) for the different preconditioners constructed using DNN, with a variable number of approximate eigenpairs retained, and with different value of the parameter $\mu$. This shows that in most cases, the preconditioner helps reach the threshold in fewer iteration than the baseline, especially for a lower value of $r$.
Furthermore, we can see that for $r=1000$ eigenpairs, the choice of $\mu$ has a limited influence on the number of iterations needed. For a larger number of eigenpairs retained (ie larger $r$), the performances are much more dependent on the the value of $\mu$, and might even reach worsen the performances (for instance $r=2000$ and $\mu$ as the minimum of the approximated eigenvalues)
The comparison of the $L_2$ norm of the residuals for the different problems is shown \cref{fig:norm_res}, which leads to similar conclusion. 
For good combinations of the parameters $r$ and $\mu$, we could reduce the number of iterations required to reach the threshold  by roughly 30\%. However, when too many eigenpairs are kept, the performances decrease.

\begin{figure}
    \centering
\begin{subfigure}[b]{0.85\textwidth}
   \includegraphics[width=1\linewidth]{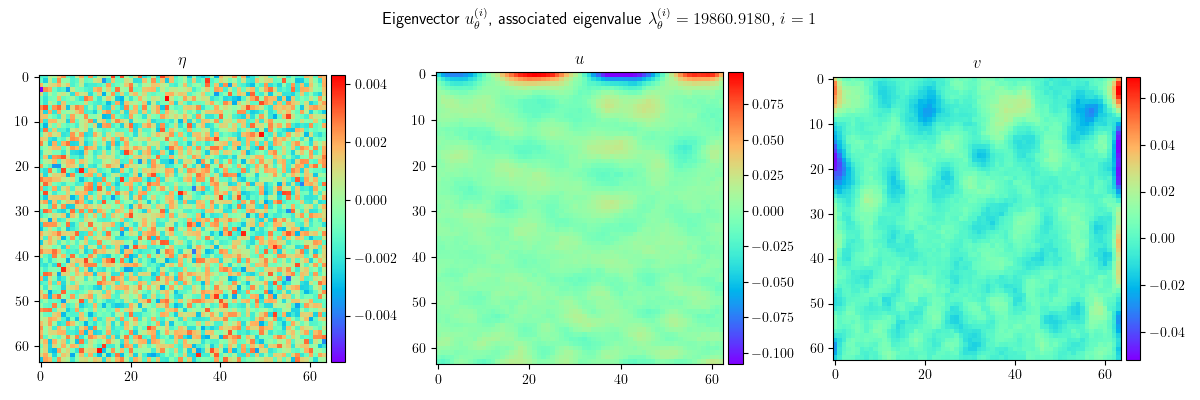}
   \caption{}
   \label{fig:_state_0} 
\end{subfigure}
\begin{subfigure}[b]{0.85\textwidth}
   \includegraphics[width=1\linewidth]{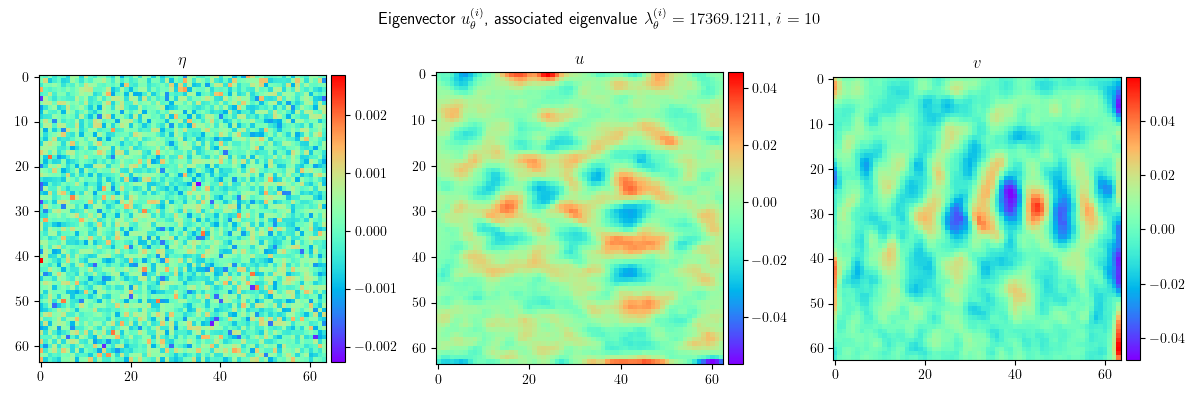}
   \caption{}
   \label{fig:_state_9} 
\end{subfigure}
\begin{subfigure}[b]{0.85\textwidth}
   \includegraphics[width=1\linewidth]{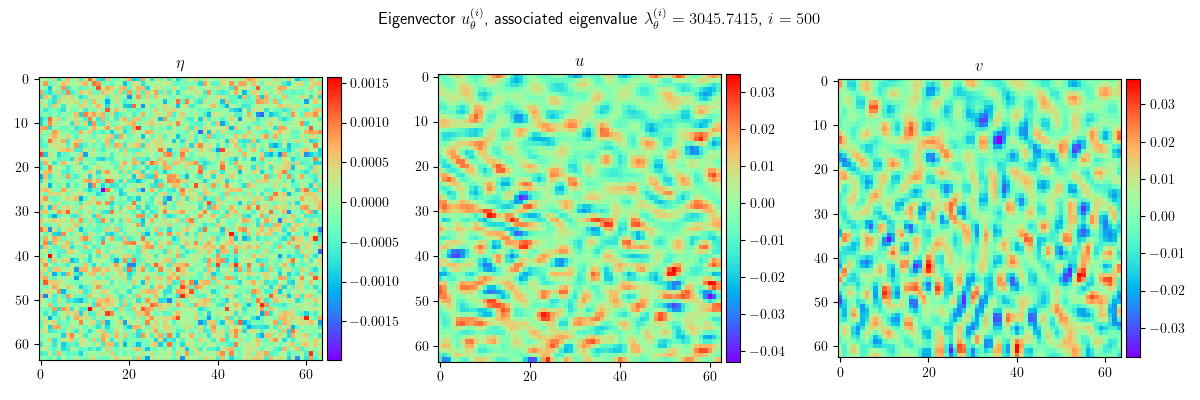}
   \caption{}
   \label{fig:_state_499} 
\end{subfigure}
\begin{subfigure}[b]{0.85\textwidth}
   \includegraphics[width=1\linewidth]{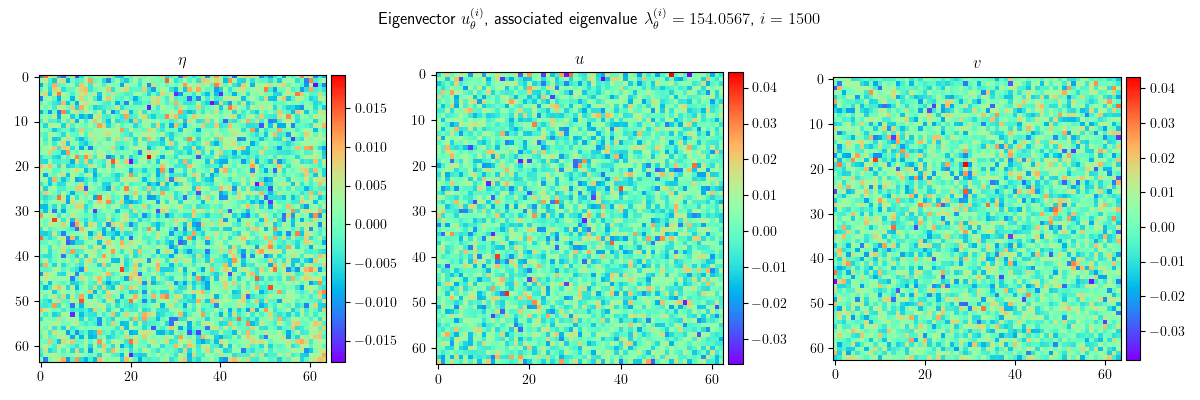}
   \caption{}
   \label{fig:_state_1499} 
\end{subfigure}
\caption{Example of estimated eigenpairs ($i=1$, $10$, $500$, $1500$) at a linearization point}
\label{fig:dnn_states}
\end{figure}

This counterintuitive result can be explained.
It is worth noting that due to the form of the reconstruction Eq.~\cref{eq:def_ml_rk1_approx}, the individual contribution of each eigenpairs gets smaller and smaller, making them more and more difficult to approximate.
We can see on \cref{fig:dnn_states} some examples of eigenvectors that the DNN outputs. The eigenvector corresponding to a rank $i=1500$ does not show any discernible pattern, in contrast with the other eigenvectors, with lower rank.
A bad estimation of the eigenpairs might worsen the quality of the preconditioner  since any error would get amplified by taking its inverse (through the negative exponent $\alpha$).
We can see this effect on the preconditioners built with $r=2000$, the whole estimated spectrum. Some of the smallest eigenvalues are not well represented by the Neural Network, and this worsen the preconditioning effect of $\precmat{\alpha}$, compared to $r=1000$ or $r=1900$, and the influence of the shift parameter $\mu$ is amplified. 
Indeed, $\mu$ helps mitigate this issue due to the approximation error of the Deep Neural Network, by forcing the resulting eigenvalues to be larger than $1$, which acts as a lower bound for the eigenvalues of the original matrix.

\section*{Conclusion and perspectives}
In this work, we focused on the problem of data-driven preconditioning of non-sparse parameterized matrices. In a Data Assimilation context, more specifically in the incremental formulation of 4D-Var, multiple resolution of high-dimensional linear systems have to be performed. For computational reasons, only a limited number of iterations of Conjugate Gradient can be done. In order to improve the rate of convergence of this iterative solving method, we propose to use Deep Neural Networks to get an approximation of the largest eigenpairs of the matrix to inverse, and then use those to precondition the linear system.

We applied this method to an academic assimilation system of moderate size. Based on the image-like structure of the state vector, we used an architecture based on U-Nets to construct a surrogate. Numerically, using this preconditioner allows for reducing the number of matrix-vector products required to reach a convergence threshold. The number of eigenpairs to use is up to the user, but a bad approximation of the eigenpairs can lead to bad performances if the parameter $\mu$ is too small. 

Compared to traditional preconditioning methods, training such a neural network can be done in a almost non-intrusive way. Once trained, this can be used as a first-level preconditioner, and thus traditional randomized methods can be applied to improve furthermore the convergence rates. 

We focused on an assimilation system where the observation operator $\mathcal{H}$ is linear and constant for all assimilation windows. In this case, the dependence on the state variable comes only from the Tangent Linear model. Because of this, the learned eigenpaires are tied to this constant observation operator (and covariance matrices of the errors). However, if those quantities were to be uniquely dependent on the state, the whole construction of the dataset and the training does not need any modification. One possible improvement of this method would be to consider a changing observation operator, and the DNN would take as input both the state of linearization, and the observation operator.
\section*{Acknowledgement}
This work has been funded within the France Relance Economic plan, and has been jointly done between Eviden and Inria.
\clearpage
\bibliographystyle{apalike}
\bibliography{bibzotero.bib}
\end{document}